\documentclass[11pt]{article}

\usepackage{amsbsy}
\usepackage{amssymb}
\usepackage{amsthm}
\usepackage{amsmath}
\usepackage{fullpage,hyperref,mathrsfs,txfonts,epsfig,graphicx,graphics}

\begin{document}

\newcommand{\bwr}{\boldsymbol{\wr}}


\newenvironment{nmath}{\begin{center}\begin{math}}{\end{math}\end{center}}

\newtheorem{thm}{Theorem}[section]
\newtheorem{lem}[thm]{Lemma}
\newtheorem{prop}[thm]{Proposition}
\newtheorem{cor}[thm]{Corollary}
\newtheorem{conj}[thm]{Conjecture}
\newtheorem{dfn}[thm]{Definition}
\newtheorem{prob}[thm]{Problem}
\newtheorem{ques}[thm]{Question}
\newtheorem{remark}[thm]{Remark}

\newcommand{\A}{\mathcal{A}}
\newcommand{\B}{\mathcal{B}}
\newcommand{\K}{\mathcal{K}}
\renewcommand{\H}{\mathcal{H}}
\renewcommand{\Pr}{\mathrm{Pr}}
\newcommand{\s}{\sigma}
\renewcommand{\P}{\mathcal{P}}
\renewcommand{\O}{\Omega}
\renewcommand{\S}{\Sigma}
\newcommand{\T}{\mathrm{T}}
\newcommand{\co}{\mathrm{co}}
\newcommand{\e}{\varepsilon}
\newcommand{\im}{\mathrm{i}}
\renewcommand{\l}{\lambda}
\newcommand{\U}{\mathcal{U}}
\newcommand{\calH}{\mathcal{H}}
\newcommand{\G}{\Gamma}
\newcommand{\g}{\gamma}
\renewcommand{\L}{\Lambda}
\newcommand{\hcf}{\mathrm{hcf}}
\newcommand{\F}{\mathcal{F}}
\renewcommand{\a}{\alpha}
\newcommand{\N}{\mathbb{N}}
\newcommand{\R}{\mathbb{R}}
\newcommand{\Z}{\mathbb{Z}}
\newcommand{\C}{\mathbb{C}}

\newcommand{\E}{\mathbb{E}}
\newcommand{\alp}{\alpha^*}

\newcommand{\bb}[1]{\mathbb{#1}}
\renewcommand{\rm}[1]{\mathrm{#1}}
\renewcommand{\cal}[1]{\mathcal{#1}}

\newcommand{\fin}{\nolinebreak\hspace{\stretch{1}}$\lhd$}

\title{The wreath product of $\Z$ with $\Z$ has Hilbert
compression exponent $\frac23$} 
\author{Tim Austin\\UCLA
\and Assaf Naor\footnote{Research supported in part by NSF grants
CCF-0635078 and DMS-0528387.}\\Courant Institute
\and Yuval Peres \footnote{Research supported in
part by NSF grant DMS-0605166.}\\Microsoft Research and UC
Berkeley
}
\date{}

\maketitle
\begin{abstract}
\noindent Let $G$ be a finitely generated group, equipped with the
word metric $d$ associated with some finite set of generators. The
 Hilbert compression exponent of $G$ is the supremum over all
$\alpha\ge 0$ such that there exists a Lipschitz mapping $f:G\to
L_2$ and a constant $c>0$ such that for all $x,y\in G$ we have
$\|f(x)-f(y)\|_2\ge cd(x,y)^\alpha.$ In~\cite{AGS06} it was shown
that the Hilbert compression exponent of the wreath product $\Z\bwr
\Z $ is at most $\frac34$, and in~\cite{NP07} was proved that this
exponent is at least $\frac23$. Here we show that $\frac23$ is the
correct value. Our proof is based on an application of K. Ball's
notion of Markov type.
\end{abstract}




\section{Introduction}

Let $G$ be a finitely generated group. Fix a finite set of
generators $S\subseteq G$, which we will always assume to be
symmetric (i.e. $S^{-1}=S$). Let $d$ be the left-invariant word
metric induced by $S$ on $G$. The {\bf Hilbert compression exponent}
of $G$, which we denote by $\alp(G)$, is the supremum over all
$\alpha\ge 0$ such that there exists a $1$-Lipschitz mapping $f:G\to
L_2$ and a constant $c>0$ such that for all $x,y\in G$ we have
$$\|f(x)-f(y)\|_2\ge cd(x,y)^\alpha.$$
Note that $\alp(G)$ does not depend on the choice of the finite set
of generators $S$, and is thus an algebraic invariant of the group
$G$. This notion was introduced by Guentner and Kaminker
in~\cite{GK04} as a natural quantitative measure of Hilbert space
embeddabililty in situations when bi-Lipschitz embeddings do not
exist (when bi-Lipschitz embeddings do exist the natural measure
would be the {\em Euclidean distortion}). More generally, the {\bf
compression function} of a $1$-Lipschitz mapping $f:G\to L_2$ is
defined as $$ \rho(t) \coloneqq \inf_{d(x,y)\ge t}
\|f(x)-f(y)\|_2.$$ The mapping $f$ is called a {\bf coarse
embedding} if $\lim_{t\to\infty}\rho(t)=\infty$. Coarse embeddings
of discrete groups have been studied extensively in recent years.
The Hilbert compression exponents of various  groups were
investigated in~\cite{GK04,AGS06,deCTesVal,Tess06,ADS06}---we refer
to these papers and the references therein for group-theoretical
motivation and applications.

Consider the wreath product $\Z\bwr \Z$, i.e. the group of all pairs
$(f,x)$, where $x\in \Z$ and $f:\Z\to \Z$ has finite support,
equipped with the group law $(f,x)(g,y)\coloneqq (z\mapsto
f(z)+g(z-x), x+y)$. In this note we prove that $\alp(\Z\bwr
\Z)=\frac23$. The problem of computing $\alp(\Z\bwr \Z)$ was raised
explicitly in~\cite{AGS06,Tess06,ADS06}. In~\cite{AGS06}
Arzhantseva, Guba and Sapir showed that $\alp(\Z\bwr \Z)\in
\left[\frac12,\frac34\right]$. In~\cite{Tess06} Tessera claimed to
improve the lower bound on $\alp(\Z\bwr \Z)$ to $\alp(\Z\bwr \Z)\ge
\frac23$, and conjectured that $\alp(\Z\bwr \Z)=\frac23$.
Unfortunately, Tessera's proof is flawed, as explained in Remark 1.4
of [12]; his method only yields the bound $\alpha^*(Z \wr Z) \geq
\frac13$. However, the inequality $\alp(\Z\bwr \Z)\ge \frac23$ is
correct, as shown by Naor and Peres in [12] using a different
method.  Here we obtain the matching upper bound $\alp(\Z\bwr \Z)\le
\frac23$. For the sake of completeness, in Remark~\ref{rem:embed}
below we also present the embeddings of Naor and Peres~\cite{NP07}
which establish the lower bound $\alp(\Z\bwr \Z)\ge \frac 23$.

Our proof of the upper bound $\alp(\Z\bwr \Z)\le \frac23$ is a
simple application of K. Ball's notion of {\bf Markov type}, a
metric invariant that has found several applications in metric
geometry in the past two
decades---see~\cite{Bal,Naor01,LMN02,BLMN05,NPSS06,MN06}. Recall
that a Markov chain $\{Z_t\}_{t=0}^\infty$ with transition
probabilities $a_{ij}\coloneqq\Pr(Z_{t+1}=j\mid Z_t=i)$ on the state
space $\{1,\ldots,n\}$ is {\em stationary\/} if
$\pi_i\coloneqq\Pr(Z_t=i)$ does not depend on $t$ and it is {\em
reversible\/} if $\pi_i\,a_{ij}=\pi_j\,a_{ji}$ for every
$i,j\in\{1,\ldots,n\}$. Given a metric space $(X,d_X)$ and $p\in
[1,\infty)$, we say that $X$ has {Markov type} $p$ if there exists a
constant $K>0$ such that for every stationary reversible Markov
chain $\{Z_t\}_{t=0}^\infty$ on $\{1,\ldots,n\}$, every mapping
$f:\{1,\ldots,n\}\to X$ and every time $t\in \mathbb N$,
\begin{eqnarray}\label{eq:defMarkov}
\E \big[ d_X(f(Z_t),f(Z_0))^p\big]\le K^p\,t\,\E\big[
d_X(f(Z_1),f(Z_0))^p\big].
\end{eqnarray}
The least such $K$ is called the Markov type $p$ constant of $X$,
and is denoted $M_p(X)$.

The fact that $L_2$ has Markov type $2$ with constant $1$, first
noted by K. Ball~\cite{Bal}, follows from a simple spectral argument
(see also inequality (8) in~\cite{NPSS06}). Since for $p\in [1,2]$
the metric space $\left(L_p,\|x-y\|_2^{p/2}\right)$ embeds
isometrically into $L_2$ (see~\cite{WW75}), it follows that $L_p$
has Markov type $p$ with constant $1$. For $p>2$ it was shown
in~\cite{NPSS06} that $L_p$ has Markov type $2$ with constants
$O\left(\sqrt{p}\right)$. We refer to~\cite{NPSS06} for a
computation of the Markov type of various additional classes of
metric spaces.

The notion of Markov type has been successfully applied to various
embedding problems of {\em finite} metric spaces. In this note we
observe that one can use this invariant in the context of infinite
amenable groups as well. In a certain sense, our argument simply
amounts to using Markov type asymptotically along neighborhoods of
F\o lner sequences.

{}For the rest of the paper, Let $G$ be an amenable group with a
fixed finite symmetric set of generators $S$ and the corresponding
left-invariant word metric $d$. Let $e$ denote the identity element
of $G$, and let $\{W_t\}_{t=0}^\infty$ be the canonical simple
random walk on the Cayley graph of $G$ determined by $S$, starting
at $e$. Our main result is:

\begin{prop}\label{prop:amenable} Assume that that there exist
$c,\delta,\beta>0$ such that for all $t\in \N$,
\begin{eqnarray}\label{eq:assumption} \Pr\left(d(W_t,e)\ge ct^\beta\right)\ge
\delta.
\end{eqnarray}
Let $(X,\, d_X)$ be a metric space with Markov type $p$, and assume
that $f:G\to X$ satisfies
\begin{eqnarray}\label{eq:compression}
\rho(d(x,y))\le d_X(f(x),f(y))\le d(x,y)
\end{eqnarray}
for all $x,y\in G$, where $\rho:\R_+\to \R_+$ is non-decreasing.
Then for all $t\in \N$,
$$
\rho\left(ct^\beta\right)\le \frac{M_p(X)}{\delta^{1/p}}t^{1/p}.
$$
In particular,
$$
\alp(G)\le \frac{1}{2\beta}.
$$
\end{prop}

As an immediate corollary we deduce that $\alp(\Z\bwr \Z)\le
\frac23$. Indeed, $\Z\bwr \Z$ is amenable (see for
example~\cite{KV83,Pat88}), and it was shown by Revelle
in~\cite{Rev03} that $\Z\bwr \Z$ has a set of generators (namely the
canonical generators $S=\{(1,0),(-1,0),(0,1),(0,-1)\}$) which
satisfies the assumption of Proposition~\ref{prop:amenable} with
$\beta=\frac34$ (see also~\cite{Ersh01} for the corresponding bound
on the expectation of $d(W_t,e)$).

\section{Proof of Proposition~\ref{prop:amenable}}

 Let
$\{F_n\}_{n=0}^\infty$ be a F\o lner sequence for $G$, i.e., for
every $\e>0$ and any finite $K\subseteq G$, we have $|F_n\triangle
(F_nK)|\le \e|F_n|$ for large enough $n$. Fix an integer $t>0$ and denote
$$
A_n\coloneqq \bigcup_{x\in F_n}B(x,t)\supseteq F_n,
$$
where $B(x,t)$ is the ball of radius $t$ centered at $x$ in the word metric determined by $S$.

For every $\e>0$, there exists $n\in \N$ such
that
\begin{eqnarray}\label{eq:size}
\e|F_n|\ge |F_n\triangle (F_nB(e,t))|=|A_n\setminus F_n| \,.
\end{eqnarray}

Let $\{Z_t\}_{t=0}^\infty$ be the delayed standard random walk
restricted to $A_n$. In other words, $Z_0$ is uniformly distributed on
$A_n$, and for all $j\ge 0$ and $x\in A_n$,
$$
\Pr\left(Z_{j+1}=x\big|Z_j=x\right)=1-\frac{|(xS)\cap A_n|}{|S|},
$$
and if $s\in S$ is such that $xs\in A_n$ then
$$
\Pr\left(Z_{j+1}=xs\big|Z_j=x\right)=\frac{1}{|S|}.
$$

It is straightforward to check that $\{Z_t\}_{t=0}^\infty$ is a
stationary reversible Markov chain. Hence, using the Markov type $p$
property of $X$, and the fact that $f$ is $1$-Lipschitz, we see that
\begin{eqnarray}\label{eq:upper}
\E\big[
d_X(f(Z_t),f(Z_0))^p\big]\stackrel{\eqref{eq:defMarkov}}{\le}
M_p(X)^pt \E\big[
d_X(f(Z_1),f(Z_0))^p\big]\stackrel{\eqref{eq:compression}}{\le}
M_p(X)^pt\E \big[d(Z_1,Z_0)^p\big]\le M_p(X)^pt.
\end{eqnarray}
Note that
\begin{equation}\label{eq:restrict}
\E \big[
d_X(f(Z_t),f(Z_0))^p\big]\stackrel{\eqref{eq:compression}}{\ge}
\E\big[ \rho(d(Z_t,Z_0))^p\big]
\ge \frac{1}{|A_n|}\sum_{x\in F_n}\E\left[
\rho\left(d(Z_t,Z_0)\right)^p\big|Z_0=x\right],
\end{equation}
since the omitted summands corresponding to $x \notin F_n$ are nonnegative.
If $x\in F_n$ then $B(x,t)\subseteq A_n$; this implies that
conditioned on the event $\{Z_0=x\}$, the random variable
$d(Z_t,Z_0)$ has the same distribution as the random variable
$d(W_t,e)$. The assumption \eqref{eq:assumption} yields that
\begin{equation} \label{eq:rw1} \E \big[\rho\left(d(W_t,e)\right)^p\big]  \ge
\rho\left(ct^\beta\right)^p\cdot\Pr\left(d(W_t,e)\ge
ct^\beta\right)\ge
 \rho\left(ct^\beta\right)^p\cdot \delta \,.
\end{equation}

In conjunction with~\eqref{eq:restrict}, this gives that
\begin{equation}\label{eq:lower}
\E \big[ d_X(f(Z_t),f(Z_0))^p\big]\ge \frac{|F_n|}{|A_n|} \cdot
\E \big[\rho\left(d(W_t,e)\right)^p\big] \stackrel{\eqref{eq:rw1}}{\ge}
\frac{|F_n|}{|A_n|}\cdot \rho\left(ct^\beta\right)^p\cdot \delta
\stackrel{\eqref{eq:size}}{\ge}
\frac{\delta}{1+\e}\cdot\rho\left(ct^\beta\right)^p.
\end{equation}
Combining~\eqref{eq:upper} and~\eqref{eq:lower}, and letting $\e\to
0$, concludes the proof of Proposition~\ref{prop:amenable}. \qed

\begin{remark}\label{rem:Q}{\em Given two groups $G$ and $H$,
the wreath product $G\bwr H$ is the group of all pairs $(f,x)$ where
$f:H\to G$ has finite support (i.e. $f(z)$ is the identity of $G$
for all but finitely many $z\in H$) and $x\in H$, equipped with the
product $ (f,x)(g,y)\coloneqq \left(z\mapsto
f(z)g(x^{-1}z),xy\right)$.
Consider
the iterated wreath products $\Z_{(k)}$, where $\Z_{(1)}=\Z$ and
$\Z_{(k+1)}\coloneqq \Z_{(k)}\bwr \Z$. In~\cite{Rev03} it is shown
that $\Z_{(k)}$ has a finite symmetric set of generators which
satisfies the assumption of Proposition~\ref{prop:amenable} with
$\beta= 1-2^{-k}$. Thus $\alp(\Z_{(k)})\le \frac{1}{2-2^{1-k}}$. In
fact, as shown in~\cite{NP07}, $\alp(\Z_{(k)})=
\frac{1}{2-2^{1-k}}$.} \fin
\end{remark}

\begin{remark}\label{rem:embed} {\em In~\cite{NP07} the lower bound $\alpha^*(\Z\bwr \Z)\ge \frac23$
is a particular case of a more general result. For the readers'
convenience we present the resulting embeddings in the case of the
group $\Z\bwr \Z$.

In what follows $\lesssim$ and $\gtrsim$ denote the corresponding
inequality up to a universal constant. Fix $\alpha\in (0,1/2)$ and
let $$ \Big\{v_g:\ g:A\to \Z\ \mathrm{finitely\  supported},\ A\in
\{\Z\cap [n,\infty)\}_{n\in \Z}\cup \{\Z\cap (-\infty,n]\}_{n\in
\Z}\Big\} $$ be disjointly supported unit vectors in $L_2(\R)$. For
$(f,k)\in \Z\bwr \Z$ define a function $\phi_\alpha(f,k): \R\to \R$
by
$$
\phi_\alpha(f,k)\coloneqq \sum_{n> k} (n-k)^\alpha \cdot
v_{f\upharpoonright_{[n,\infty)}}+\sum_{n< k} (k-n)^\alpha\cdot
v_{f\upharpoonright_{(-\infty,n]}}.
$$
Observe that $\phi_\alpha(f,k)-\phi_\alpha(0,0)\in L_2(\R)$. Indeed,
if $f$ is supported on $[-m,m]$ then
$$
\|\phi_\alpha(f,k)-\phi_\alpha(0,0)\|_2^2\lesssim
m\left(m^{2\alpha}+|k|^{2\alpha}\right)+\sum_{n\in \Z}
\left(|n|^{\alpha}-|n-k|^{\alpha}\right)^2\lesssim
m\left(m^{2\alpha}+|k|^{2\alpha}\right)+\sum_{j=1}^\infty\frac{k^2}{j^{2(1-\alpha)}}
<\infty.
$$
We can therefore define $F_\alpha:\Z\bwr \Z\to \R\oplus
\ell_2(\Z)\oplus L_2(\R)$ by $$ F_\alpha(f,k)\coloneqq k\oplus
f\oplus \big(\phi_\alpha(f,k)-\phi_\alpha(0,0)\big).$$  We claim
that for every $(f,k)\in \Z\bwr \Z$ we have
\begin{eqnarray}\label{eq:biLip2/3}
d_{\Z\bwr
\Z}\big((f,k),(0,0)\big)^{\frac{2\alpha+1}{2\alpha+2}}\lesssim
\|F_\alpha(f,k)\|_2\lesssim \frac{1}{\sqrt{1-2\alpha}}\cdot
d_{\Z\bwr \Z}\big((f,k),(0,0)\big),
\end{eqnarray}
Since the metric $\|F_\alpha(f_1,k_1)-F_\alpha(f_2,k_2)\|_2$ is
$\Z\bwr\Z$-invariant, and $F_\alpha(0,0)=0$, the inequalities
in~\eqref{eq:biLip2/3} imply that $\Z\bwr \Z$ has Hilbert
compression exponent at least $\frac{2\alpha+1}{2\alpha+2}$. Letting
$\alpha\uparrow \frac12$ shows that $\alpha^*(\Z\bwr \Z)\ge
\frac23$.

It suffices to check the upper bound in~\eqref{eq:biLip2/3} (i.e.
the Lipschitz condition for $F_\alpha$) when $(f,k)$ is one of the
generators of $\Z\bwr \Z$, i.e. $(f,k)=(0,1)$ or
$(f,k)=(\delta_0,0)$. Observe that $ \|F_\alpha(\delta_0,0)\|_2=1$
and
$$
\|F_\alpha(0,1)\|_2^2\lesssim
\sum_{n=1}^\infty\left(n^\alpha-(n-1)^\alpha\right)^2\lesssim
\frac{1}{1-2\alpha},
$$
implying the upper bound in~\eqref{eq:biLip2/3}. To prove the lower
bound in~\eqref{eq:biLip2/3} assume that $m\in \N$ is the minimal
integer such that $f$ is supported on $[k-m,k+m]$. Then,
\begin{multline*}
\|F_\alpha(f,k)\|^2_2\gtrsim k^2+\sum_{j=k-m}^{k+m}
f(j)^2+\sum_{\ell=1}^m \ell^{2\alpha}\gtrsim
k^2+\frac{1}{m}\left(\sum_{j\in \Z}
|f(j)|\right)^2+m^{2\alpha+1}\\\gtrsim \left(k+m+\sum_{j\in \Z}
|f(j)|\right)^{\frac{4\alpha+2}{2\alpha+2}}\gtrsim d_{\Z\bwr
\Z}\big((f,k),(0,0)\big)^{\frac{4\alpha+2}{2\alpha+2}},
\end{multline*}
where the penultimate inequality follows by considering the cases
$\|f\|_1\ge m^{\alpha+1}$ and $\|f\|_1\le m^{\alpha+1}$ separately.}
\fin
\end{remark}

\bigskip

\bigskip


\bibliographystyle{abbrv}
\bibliography{wreath}

\end{document}